\documentclass[11 pt,reqno]{amsart}
\usepackage{amssymb}
\usepackage{amsmath}
\input xypic

\newcommand{\can}{\mathrm{can}}
\newcommand{\pos}{\operatorname{pos}}
\def\iso{\buildrel \sim\over\to}

\newcommand {\Ad}{\operatorname{Ad}}
\newcommand {\Alt}{\operatorname{Alt}}
\newcommand {\Aut}{\operatorname{Aut}}

\newcommand {\Ext}{\operatorname{Ext}}
\newcommand {\Hom}{\operatorname{Hom}}
\renewcommand {\Im}{\operatorname{Im}}

\newcommand {\Ker}{\operatorname{Ker}}

\newcommand {\qca}{quasi--Coxeter algebra }
\newcommand {\qcas}{quasi--Coxeter algebras }
\newcommand {\rhs}{right--hand side }
\newcommand {\qis}{quasi--isomorphism }
\newcommand {\qiss}{quasi--isomorphisms }
\renewcommand {\inf}{_{inf}}

\newcommand {\sfA}{\mathsf{A}}
\newcommand {\sfB}{\mathsf{B}}
\newcommand {\sfD}{\mathsf{D}}
\newcommand {\sfE}{\mathsf{E}}
\newcommand {\sfF}{\mathsf{F}}
\newcommand {\sfG}{\mathsf{G}}
\newcommand {\sfH}{\mathsf{H}}
\newcommand {\sfI}{\mathsf{I}}
\newcommand {\sfX}{\mathsf{X}}
\newcommand {\CC}{\mathcal C}
\newcommand {\CF}{\mathcal F}
\newcommand {\CO}{\mathcal O}
\newcommand {\CP}{\mathcal P}
\newcommand {\GS}{\mathfrak S}
\newcommand {\IC}{\mathbb{C}}
\newcommand {\IN}{\mathbb{N}}                          
\newcommand {\IZ}{\mathbb{Z}}
\newcommand {\veps}{\varepsilon}
\newcommand {\cupone}{\cup\{1\}}
\newcommand {\tp}{_{2,+}}
\newcommand {\II}{{\scriptscriptstyle{\operatorname{II}}}}

\newcommand {\ths}{\thinspace }
\newcommand {\wt}{\widetilde}

\newcommand {\aalpha}{\underline{\alpha}}
\newcommand {\bbeta}{\underline{\beta}}
\newcommand {\overlineambda}{{\overline{\lambda}}}
\newcommand {\aand}{\qquad\text{and}\qquad}
\newcommand {\con}[3]{{\complement}^{#1\setminus #2}_{#3}}
\newcommand {\odots}[1]{#1\cdots #1}

\newcommand {\rot}[1]{\alpha_{#1}}

\newcommand {\alg}[1]{{k#1}}
\renewcommand {\vec}[1]{{k\left(#1\right)}}

\newtheorem {theorem}{Theorem}[section]
\newtheorem {proposition}[theorem]{Proposition}

\newtheorem {lemma}[theorem]{Lemma}

\newtheorem {definition}[theorem]{Definition}

\newcommand {\remark}
{\addtocounter{theorem}{1}{\noindent\bf Remark \thetheorem.}\ }

\newcommand {\fd}{finite--dimensional }
\newcommand {\g}{\mathfrak g}
\renewcommand {\H}{\mathcal H}

\newcommand {\Gr}{\mathsf{Gr}}

\newcommand {\R}{R}
\newcommand {\pprime}{^\prime}
\newcommand {\fml}{[\negthinspace[\hbar]\negthinspace]}

\begin{document}

\title{The Dynkin diagram cohomology of finite Coxeter groups}
\author[R. Rouquier]{Rapha\"{e}l Rouquier}
\address{
R.R.: Mathematical Institute, University of Oxford\\
24-29 St Giles'\\
Oxford, OX1 3LB, UK}
\email{rouquier@maths.ox.ac.uk}
\author[V. Toledano Laredo]{Valerio Toledano Laredo}
\thanks{V.T.L. is supported in part by NSF grants DMS--0707212 and DMS--0635607}
\address{
V.T.L.: Department of Mathematics\\
Northeastern University\\
567 Lake Hall\\
360 Huntington Avenue\\
Boston\\
MA 02115
\newline\newline
School of Mathematics\\
Institute for Advanced Study\\
Princeton, NJ 08540}
\email{V.ToledanoLaredo@neu.edu}
\begin{abstract}
Let $D$ be a connected graph. The Dynkin complex $CD(A)$
of a $D$--algebra $A$ was introduced by the second author in
\cite{TL2} to control the deformations of \qca structures on $A$.
In the present paper, we study the cohomology of this complex 
when $A$ is the group algebra of a Coxeter group $W$ and
$D$ is the Dynkin diagram of $W$. We compute this cohomology
when $W$ is finite and prove in particular the rigidity of \qca
structures on $kW$. For an arbitrary $W$, we compute the top
cohomology group and obtain a number of additional partial
results when $W$ is affine. Our computations are carried out
by filtering the Dynkin complex by the number of vertices of
subgraphs of $D$. The  corresponding graded complex turns
out to be dual to the sum of the Coxeter complexes of all standard,
irreducible parabolic subgroups of $W$.
\end{abstract}
\maketitle


\section{Introduction}

Let $\g$ be a complex, semisimple Lie algebra and $D$ the
corresponding Dynkin diagram.
The notion of quasi--Coxeter algebra of type $D$ was introduced
in \cite{TL2} to put the monodromy of the Casimir connection of $
\g$ \cite{MTL,TL} and the quantum Weyl group representations
arising from the quantum group $U_q\g$ \cite{Lu} on an equal
footing, and allow for their comparison via a suitable deformation
complex.

Roughly speaking, a \qca $A$ of type $D$ is an algebra which
carries representations of the generalised braid group $B_D$
corresponding to $D$ on its \fd modules. The deformation theory
of such an algebra is controlled by a complex $CD(A)$ concentrated
in degrees $0\leq p\leq |D|$ called the {\it Dynkin complex} of $A$.

If $W$ is the Weyl group of $\g$, the complex group algebra
$A=\IC W$ possesses \qca structures accounting respectively
for the representations of $B_D$ coming from the Hecke algebra
$\H(W)$ of $W$ and the monodromy of Cherednik's KZ
connection \cite{Ch} (see \cite[\S 4]{TL2} for details).

While in this case these representations may easily be shown
to be equivalent, 
this raises nonetheless the question of computing the Dynkin
diagram cohomology of $\IC W$ and more generally of the
group algebra $k W$ of a Coxeter group over an arbitrary
ground ring $k$, when the underlying $D$--algebra structure
arises from the standard parabolic subgroups of $W$.

In the present paper, we answer this question when $W$ is
finite and $k$ is a field of characteristic 0. Along the way, we
also determine the top cohomology groups for arbitrary Coxeter
groups and obtain partial results for affine ones.

We carry out our computations by filtering $CD(kW)$ by the
number of vertices of subgraphs of $D$. Interestingly, and
crucially for us, the associated graded complex $\Gr(CD(kW))$
turns out to be dual to the sum of the Coxeter complexes of all
standard, irreducible, parabolic subgroups of $W$. This greatly
simplifies the computation of $H^*(CD(kW))$ since the Coxeter
complex is acyclic for affine Coxeter groups and has cohomology
in one degree only for finite ones.

We turn now to a detailed description of the paper.

In Section \ref{ss:Dalg}, we review the definition of the Dynkin
complex of a $D$--algebra and define its canonical filtration.

We then consider in Section \ref{se:filtration} the case of the
group algebra of a Coxeter group. In this case, the
Dynkin complex has a direct sum decomposition parametrised
by the conjugacy classes of $W$. Further, the graded complex
$\Gr(CD(kW))$ determined by the filtration of $CD(kW)$
decomposes as a sum parametrised by connected subgraphs
$B$ of $D$. The summand associated to $B$ is given by the
morphisms from the Coxeter complex of $W_B$ to $kW_B$
endowed with the adjoint action of $W_B$.

In Section \ref{sec:cohCDB}, we use the known description of
the cohomology of the Coxeter complex for finite and affine
Coxeter groups to compute the cohomology of $\Gr(CD(kW))$.

In Section \ref{se:E2}, we consider finite Coxeter groups.
We compute the induced differential on $H^*(\Gr(CD(kW)))$
and show that the resulting complex is quasi--isomorphic
to $CD$.

We apply this result in Section \ref{se:qca} to prove the rigidity
of \qca structures on $kW$ when $W$ is finite.

The main result of Section \ref{sec:top} is the construction of a
basis of the top Dynkin cohomology of the group algebra of an
arbitrary Coxeter group: it is parametrised by cuspidal conjugacy
classes such that the centraliser of an element of the class is in
the kernel of the sign character.

Section \ref{sec:finite} is devoted to the determination of the Dynkin
cohomology of finite Coxeter groups. We proceed case by case,
by explicit computation for $W$ classical or of type $\sfG_2$ and
$\sfF_4$, and by using the program GAP for the remaining exceptional
groups. In types $\sfA$ and $\sfB$, we provide a very simple formula
for the corresponding generating series which, in type $\sfA$ turns
out to be the product of a bosonic and fermionic partition function.
We show moreover that, for $W$ classical, the Dynkin cohomology
spaces stabilise with the rank of $W$.

In the final Section \ref{sec:affine}, we describe the part of the
Dynkin cohomology of an affine Weyl group corresponding to
conjugacy classes of elements of infinite order.

\section{$D$-algebras and the Dynkin complex}
\label{ss:Dalg}

This section reviews the definition of $D$-algebras and of the Dynkin complex.
With the exception of \S \ref{ss:filtration}, the material is borrowed from \cite{TL2}.

\subsection{$D$-algebras \cite[\S 3]{TL2}}

Let $D$ be a connected {\it diagram}, that is a non-empty undirected graph
with no multiple edges or loops. We denote the set of vertices of $D$ by $V(D)$
and set $|D|=|V(D)|$. By a {\it subdiagram} $B\subset D$ we shall mean a
non-empty
full subgraph of $D$, that is a graph consisting of a subset $V(B)$ of vertices of
$D$, together with all edges of $D$ joining any two elements of $V(B)$. We
will often abusively identify such a $B$ with its set of vertices and write $\alpha
\in B$ to mean $\alpha\in V(B)$. Two subdiagrams $B_1,B_2\subseteq D$ are
{\it orthogonal} if they have no vertices in common and no two vertices
$\alpha_1\in B_1,\alpha_2\in B_2$ are joined by an edge in $D$.\\

Let $k$ be a fixed commutative ring with unit. By an algebra we shall mean a
unital, associative $k$-algebra. All algebra homomorphisms are assumed to
be unital. Recall the following

\begin{definition}
A $D$-algebra is an algebra $A$ endowed with subalgebras $A_B$
labelled by the connected subdiagrams $B$ of $D$ such that the
following holds
\begin{itemize}
\item $A_B\subseteq A_{B'}$ whenever $B\subseteq B'$.
\item $A_B$ and $A_{B'}$ commute whenever $B$ and $B'$ are orthogonal.
\end{itemize}
\end{definition}

If $B_1,B_2\subseteq D$ are subdiagrams with $B_1$ connected, we denote
by $A_{B_1}^{B_2}$ the centraliser in $A_{B_1}$ of the subalgebras $A_{B_
2^i}$, where $B_2^i$ runs over the connected components of $B_2$.

\subsection{The Dynkin complex \cite[\S 5]{TL2}}

For any $0\leq p\leq n=|D|$, set
$$C^{p}(A)=
\bigoplus_{\substack{\aalpha\subseteq B\subseteq D,\\|\aalpha|=p}}
A_B^{B\setminus\aalpha}$$
where the sum ranges over all connected subdiagrams $B$ of $D$ and ordered
subsets $\aalpha=\{\alpha_1,\ldots,\alpha_p\}\subseteq V(B)$ of cardinality $p$
and, by convention, $A_B^\emptyset=A_B$. We denote the component of $a\in
C^p(A)$ along $A_B^{B\setminus\aalpha}$ by $a_{(B;\aalpha)}$.

\begin{definition} The group of Dynkin $p$-cochains on $A$ is the subspace
$CD^{p}(A)\subset C^{p}(A)$ of elements $a$ such that $$a_{(B;\sigma\aalpha)}
=(-1)^{\sigma}a_{(B;\aalpha)}$$ where, for any $\sigma\in\GS_{p}$, $\sigma
\{\alpha_1,\ldots,\alpha_p\}=\{\alpha_{\sigma(1)},\ldots,\alpha_{\sigma(p)}\}$
\footnote{In \cite{TL2} Dynkin chains are defined with values in any $D$--bimodule
$M$ over $A$. We shall only need to consider $M=A$ in this paper and therefore
restrict attention to this case.}
\end{definition}

Note that
$$CD^0(A)=\bigoplus_{B\subseteq D}Z(A_B)
\aand
CD^{n}(A)\simeq A_D$$
For $1\leq p\leq n-1$, define a map $d_D^{p}:C^{p}(A)\rightarrow
C^{p+1}(A)$ by
\begin{equation}\label{eq:Dynkin differential}
\left(d_D a\right)_{(B;\aalpha)}=
\sum_{i=1}^{p+1}(-1)^{i-1}\left(
a_{(B;\aalpha\setminus\rot{i})}-
a_{(\con{B}{\rot{i}}{\aalpha\setminus\rot{i}};
\aalpha\setminus\rot{i})}
\right)
\end{equation}
where $\aalpha=\{\rot{1},\ldots,\rot{p+1}\}$,
$\con{B}{\rot{i}}{\aalpha\setminus\rot{i}}$ is
the connected component of $B\setminus\rot{i}$
containing $\aalpha\setminus\rot{i}$ if one such
exists and the empty set otherwise, and we set $
a_{(\emptyset;-)}=0$. For $p=0$, define $d_D^0:
C^{0}(A)\rightarrow C^{1}(A)$ by
$$d_D^0 a_{(B;\alpha_i)}=a_B-a_{B\setminus\rot{i}}$$
where $a_{B\setminus\alpha_i}$ is the sum of $a_{B_2}$
with $B_2$ ranging over the connected components of
$B\setminus\rot{i}$. Finally, set $d_D^n=0$. The map $d_D$ leaves
$CD(A)$ invariant and satisfies $d_D^2=0$. The cohomology
$HD(A)$ of $CD(A)$ with respect to $d_D$ is called the
{\it Dynkin diagram} cohomology of $A$. 

\subsection{Restriction \cite[\S 5.2]{TL2}}\label{ss:restriction}

Let $D'$ be a connected subgraph of $D$.
We have a morphism of complexes
\begin{align*}
\mathrm{Res}_{D'}^D:CD(A)&\to CD(A_{D'}) \\
A_B^{B\setminus\underline{\alpha}}\ni a&\mapsto
\begin{cases}
a & \textrm{ if }B\subset D' \\
0 & \textrm{ otherwise.}
\end{cases}
\end{align*}

\subsection{The canonical filtration on $CD(A)$}\label{ss:filtration}

Endow each chain group $C^p(A)$ with the $\IN$-grading given by
$$C^p_q(A)=
\bigoplus_
{\substack{\aalpha\subset B\subset D\\|\aalpha|=p,\medspace|B|=q}}
A_B^{B\setminus\aalpha}$$
where $p\leq q\leq n$, and set $CD^p_q(A)=CD^p(A)\cap C^p_q(A)$.
Since $|\con{B}{\rot{i}}{\aalpha\setminus\rot{i}}|< |B|$, the Dynkin
differential $d_D$ maps $CD^p_q(A)$ to
$$CD^{p+1}_{\geq q}(A)=\bigoplus_{r=q}^n CD^{p+1}_r(A)$$
This gives a decreasing $\IN$--filtration on the Dynkin complex of $A$.
The $\sfE_2$--term of the corresponding spectral sequence is the
cohomology of $CD(A)$ with respect to the differential
\begin{equation}\label{eq:d0}
\left(d_D^0 a\right)_{(B;\aalpha)}=
\sum_{i=1}^{p+1}(-1)^{i-1}
a_{(B;\aalpha\setminus\rot{i})}
\end{equation}

\section{The Dynkin complex of a Coxeter group}
\label{se:filtration}

\subsection{Description}

Let $W$ be an irreducible Coxeter group with system of generators
$S=\{s_i\}_{i\in I}$ and let $D$ be the Coxeter graph of $(W,S)$. For
any subgraph $B\subseteq D$ with vertex set $I_B\subseteq I$, let
$W_B\subseteq W$ be the standard parabolic subgroup generated
by $s_i$, $i\in I_B$.\\

Regard the group algebra $A=\alg{W}$ as a $D$-algebra by setting
$A_B=\alg{W_B}$. By choosing a total order on the vertices of $D$,
we can identify the corresponding Dynkin complex with
$$CD^p=
\bigoplus_{\substack{\aalpha\subset B\subset D\\|\aalpha|=p}}
\alg{W_B}^{W_{B\setminus\aalpha}}$$
where $\alpha$ now ranges over the unordered subsets of $V(B)$
and the $W_{B\setminus\aalpha}$--fixed points in $\alg{W_B}$ are
taken with respect to the diagonal (adjoint) action.
The Dynkin differential on $\alg{W_B}^{W_{B\setminus\aalpha}}$
is the map
$$\alg{W_B}^{W_{B\setminus\aalpha}}\longrightarrow
\bigoplus_{\beta\in B\setminus\aalpha}
\alg{W_B}^{W_{B\setminus(\aalpha\cup\{\beta\})}}
\oplus
\bigoplus_{\substack{B\subset B'\subset D,\\\beta\in B'\setminus\aalpha\,:\\
B=\complement^{B'\setminus\{\beta\}}_{\aalpha}}}
\alg{W_{B'}}^{W_{B'\setminus(\aalpha\cup\{\beta\})}}$$
given by
$$\left(\sum_\beta (-1)^{\mathrm{pos}(\beta,
\aalpha\cup\{\beta\})-1}\can,
\sum_{B',\beta} (-1)^{\mathrm{pos}(\beta,
\aalpha\cup\{\beta\})}\can\right)$$
where $\can$ is the canonical inclusion map and for any $\beta\in\bbeta
\subset V(D)$, $\pos(\beta,\bbeta)\in\{1,\ldots,|\bbeta|\}$ is the position
of $\beta$ relative to the total order on $\bbeta$.

\subsection{Decomposition by conjugacy classes}\label{ss:conjugacy}

Let $\CC$ be the set of conjugacy classes of $W$. For any $c\in\CC$,
set 
$$CD^p_{c}=
\bigoplus_{\substack{\aalpha\subset B\subset D\\|\aalpha|=p}}
\vec{W_B\cap c}^{W_{B\setminus\aalpha}}$$
where, for any set $X$, $\vec{X}$ is the vector space with basis $X$.
Then $CD_{c}=\bigoplus_p CD^p_{c}$ is a subcomplex of $CD$ and
$$CD=\bigoplus_{c\in \CC}CD_{c}.$$
We denote the cohomology of the corresponding complex by $HD
_c(\alg{W})$.

\subsection{Filtration}

When filtered as in \S \ref{ss:filtration}, the associated graded complex
$\Gr(CD)$ (= the $\sfE_1$-term of the spectral sequence)
is the sum over all connected subdiagrams $B\subseteq D$ of the
subcomplexes $CD_B$ given by
\begin{equation}
\label{eq:adjunction}
\begin{split}
CD^p_B &=
\bigoplus_{\aalpha\subset B,\medspace |\aalpha|=p}
\alg{W_B}^{W_{B\setminus\aalpha}}\\
&\simeq
\bigoplus_{\aalpha\subset B,\medspace |\aalpha|=p}
\Hom_{W_B}(\vec{W_B/W_{B\setminus\aalpha}},\alg{W_B})
\end{split}
\end{equation}

Recall that the Coxeter complex $CC^B$ of $W_B$ is the
(homology) complex
$$CC^B_p=
\bigoplus_{\aalpha\subset B,\medspace |\aalpha|=p}
\vec{W_B/W_{B\setminus\aalpha}}$$
with differential given by
$$\partial_C\left(wW_{B\setminus\aalpha}\right)=
\sum_{i=1}^p(-1)^{i-1}\medspace
wW_{(B\setminus\aalpha)\cup\{\alpha_i\}}$$
where $\aalpha=\{\alpha_1,\ldots,\alpha_p\}$ with $\alpha_1<\cdots
<\alpha_p$.
The following immediate result identifies the $B$--component
$CD_B$ of $\Gr(CD)$ with the dual of the Coxeter
complex of $W_B$ with values in $\alg{W_B}$.

\begin{proposition}\label{pr:Dynkin=Coxeter}
The isomorphism \emph{(\ref{eq:adjunction})}
induces an isomorphism of complexes
$$CD_B\simeq
\Hom_{W_B}(CC^B,\alg{W_B}).$$
\end{proposition}


\section{The cohomology of $CD_B$}
\label{sec:cohCDB}

Assume henceforth that $k$ is a field of characteristic 0.

\subsection{Finite Coxeter groups}\label{ss:finite CDB}

Assume in this paragraph that $W_B$ is finite. Let $S_B$ be the
unit sphere in the Euclidean reflection representation of $W_B$,
and cellulate $S_B$ by its intersections with the chambers of $
W_B$. Then, the Coxeter complex $CC^B$ is the cellular homology
complex of $S_B$, reduced and shifted by one \cite{Hu}. Thus,
$H_p(CC^B)$ is zero if $p<|B|$ and the sign representation
$\veps$ of $W_B$ otherwise, so that, by Proposition \ref{pr:Dynkin=Coxeter},
we have
\begin{equation}\label{eq:sphere}
H^p(CD_B)\simeq
\left\{\begin{array}{cl}
\alg{W_B}^\veps	&\text{if $p=|B|$}\\
0		&\text{otherwise.}
\end{array}\right.
\end{equation}
where $\alg{W_B}^\veps\subset\alg{W_B}$ is the subspace
transforming like the sign representation $\veps$ of $W_B$.

Consider $\alg{W_B}^\veps[-|B|]$, a complex concentrated in
degree $|B|$. We have a morphism of complexes $i_B:{\alg{W_B}}^
\veps[-|B|]\to CD_B$ given by the inclusion
$$\alg{W_B}^\veps\hookrightarrow
\alg{W_B}=
\alg{W_B}^{W_{B\setminus B}}$$
in degree $|B|$. Let $\Alt_B:\alg{W_B}\rightarrow \alg{W_B}^
\veps$ be the projection given by 
$$\Alt_B(f)=\frac{1}{|W_B|}\sum_{w\in W_B}\veps(w) wfw^{-1}.$$
Since $\Alt_B$ is zero on
$\sum_{\alpha\in B}\alg{W_B}^{W_{B\setminus\alpha}}$,
it defines a morphism of complexes
$\rho_B:CD_B\to \alg{W_B}^\veps[-|B|]$. Summarising, we
have the following Proposition.

\begin{proposition}\label{pr:quasiisoCoxeter}
If $W_B$ is finite, the maps $i_B:\alg{W_B}^\veps[-|B|]\to
CD_B$ and $\rho_B:CD_B\to \alg{W_B}^\veps[-|B|]$ are
\qiss such that $\rho_B\circ i_B=\mathrm{id}$.
\end{proposition}

\subsection{Affine Coxeter groups}
\label{se:affinecox}

Assume now that $W_B$ is an affine Coxeter group, and let $E_B$
be the Euclidean space of dimension $|B|-1$ cellulated by the alcoves
of $W_B$. The Coxeter complex $CC^B$ is the cellular homology
complex of $E_B$, reduced and shifted by one \cite{Hu}, and is
therefore acyclic. Its terms of positive degree are induced from
the trivial representation of finite (parabolic) subgroups of $W_B$
and are therefore projective. Thus, $CC^B$ is a projective resolution
of the trivial $W_B$--module. It follows from this, and Proposition
\ref{pr:Dynkin=Coxeter}, that
\begin{equation}\label{eq:affine Coxeter}
H^p(CD_B)\simeq
\left\{\begin{array}{cl}
\Ext^{p-1}_{W_B}(k,\alg{W_B}) &\text{if $p\geq 2$}\\
0                      &\text{otherwise.}
\end{array}\right.
\end{equation}
where $\alg{W_B}$ is endowed with the adjoint action of $W_B$.

\section{The $\sfE_2$--term of the Dynkin complex of a Coxeter group}\label{se:E2}

\subsection{Finite Coxeter groups}
\label{ss:finiteHC}

We assume in this section that the Coxeter group $W$ is finite.\\

Consider the complex
\begin{equation}\label{eq:HC}
HC=\bigoplus_{B\subseteq D}\alg{W_B}^\veps
\end{equation}
where $\alg{W_B}^{\veps}$ appears in degree $|B|$, endowed
with the differential $d^\#$ given by
\begin{equation}\label{eq:d flat}
d^\#=\sum_{B'}(-1)^{(B';B)}\cdot\Alt_{B'}:
\alg{W_B}^\veps\longrightarrow
\bigoplus_{\substack{B\subset B'\\ |B'|=|B|+1}}
\alg{W_{B'}}^\veps
\end{equation}
where, given $B'\subseteq D$ with ordered set of vertices $\alpha
_1\odots{<}\alpha_q$ and $B\subset B'$ such that $B'\setminus
B=\{\alpha_i\}$, we set $(B';B)=i$.

Note that this complex is concentrated in degrees $\geq 2$ since,
for $|B|=1$, we have $\alg{W_B}^\veps=\alg{\GS_2}^\veps=0$.

\smallskip
Consider the application $\phi:CD\to HC$ given by
\begin{equation}\label{eq:alt quiso}
\alg{W_B}^{W_{B\setminus\aalpha}}\ni a\mapsto\Alt_B(a).
\end{equation}

\begin{proposition}
\label{pr:CDtoHC}
The application $\phi$ is a \qis of complexes.
\end{proposition}

\begin{proof}
Let $a\in \alg{W_B}^{W_{B\setminus\aalpha}}$. We have
$\phi d_D(a)=0=\phi(a)$ if $\aalpha{\not=}B$.
On the other hand, if $\aalpha=B$, we have
$$\phi d_D(a)=\sum_{\substack{B\subset B'\subset D\\ |B'|=|B|+1}}
(-1)^{\mathrm{pos}(B'\setminus B,B')}\Alt(a)=d_D\phi(a).$$
so that $\phi$ is compatible with the differential.

Consider the filtration on $HC$ given by
$(HC_{\ge q})^p=HC^p$ for $p\ge q$ and $(HC_{\ge q})^p=0$ for $p<q$. The
morphism $\phi$ is a morphism of filtered complexes. Via the canonical
isomorphisms of \S \ref{se:filtration},
the induced morphism $\bar{\phi}_p:CD_{\ge p}/CD_{\ge p+1}\to
HC_{\ge p}/HC_{\ge p+1}$ becomes
the sum over connected subdiagrams $B$ of $D$ of cardinality $p$ of
the morphisms $\rho_B$ of \S \ref{ss:finite CDB}. It follows that $\bar
{\phi}_p$ is a \qis by Proposition \ref{pr:quasiisoCoxeter}. The Proposition
follows.
\end{proof}

Let $D'$ be a connected subgraph of $D$ and $HC'$ the corresponding
complex. Via the isomorphisms $\phi$ above, the restriction map of \S
\ref{ss:restriction} becomes
\begin{align*}
\mathrm{Res}_{D'}^D:HC&\to HC'\\
kW_B^\veps\ni a&\mapsto
\begin{cases}
a & \textrm{ if }B\subset D' \\
0 & \textrm{ otherwise.}
\end{cases}
\end{align*}

\remark It seems an interesting problem to determine whether the complex
\eqref{eq:HC}--\eqref{eq:d flat} is the cellular cochain of a $CW$--complex
or the Morse complex of a smooth manifold naturally associated to $W$.

\subsection{Finite part}
\label{se:finitepart}

Assume now that $W$ is an arbitrary Coxeter group. We proceed as in 
\S \ref{ss:finiteHC} for the subspace of $CD$ corresponding to subdiagrams
$B\subseteq D$ such that $W_B$ is finite.

Let $CD_{\mathrm{inf}}\subset CD$ be the subcomplex
$$CD_{\mathrm{inf}}=\bigoplus_{\substack{B\subseteq D:\\|W_B|=\infty}}CD_B$$

Let 
$$HC_f=\bigoplus_{\substack{B\subset D:\\|W_B|<\infty}}\alg{W_B}^\veps$$
where $B$ runs over connected subdiagrams of $D$ which are Dynkin.
Here, $\alg{W_B}^{\veps}$ is in degree $|B|$, endowed with
the differential $d^\#$ given by
$$\sum_{B'}(-1)^{(B';B)}\cdot\Alt_B:
\alg{W_B}^\veps
\longrightarrow
\bigoplus_{\substack{B\subset B'\\ |B'|=|B|+1}}
\alg{W_{B'}}^\veps$$

Consider the application
\begin{align*}
\phi:CD&\to HC_f\\
\alg{W_B}^{W_{B\setminus\aalpha}}\ni a&\mapsto 
\begin{cases}
\Alt_B(a) & \text{ if }\aalpha=B \text{ and $|W_B|<\infty$}\\
0 & \text{ otherwise.}
\end{cases}
\end{align*}
As in Proposition \ref{pr:CDtoHC}, one checks that $\phi$ is 
a morphism of complexes.

\begin{proposition}
There is a distinguished triangle
$$CD_{\mathrm{inf}}\xrightarrow{\can} CD\xrightarrow{\phi} HC_f\rightsquigarrow$$
\end{proposition}

\begin{proof}
One shows as in Proposition \ref{pr:CDtoHC} that the map $CD/CD\inf
\to HC_f$ induced by $\phi$ is a quasi--isomorphism.
\end{proof}

\section{Rigidity of \qca structures on $kW$}\label{se:qca}

We apply below the results of Section \ref{se:E2} to show that \qca
structures on $kW$ are rigid if $W$ is finite. We begin by reviewing
the definition of \qcas and their deformations.

\subsection{Quasi--Coxeter algebras \cite[\S 3]{TL2}} 

Let $\R$ be a commutative ring with unit. Recall that a {\it \qca} structure
on $\R W$ is given by endowing it with the following data
\begin{itemize}
\item {\bf Local monodromies:} %
for each $i\in I$, 
an invertible element $S_i\in \R W_{\alpha_i}\cong\R\IZ_2$.
\item {\bf Elementary associators:} %
for each connected subdiagram $B\subseteq D$ and vertices $\alpha_i\neq
\alpha_j\in B$, an invertible element $\Phi_{(B;\alpha_i,\alpha_j)}\in \R W$.
\end{itemize}
satisfying the following axioms
\begin{itemize}
\item {\bf Orientation:}%
$$\Phi_{(B;\alpha_j,\alpha_i)}=\Phi_{(B;\alpha_i,\alpha_j)}^{-1}$$
\item {\bf Support:}%
$$\Phi_{(B;\alpha_i,\alpha_j)}\in {\R W_B}^{B\setminus\{\alpha_i,\alpha_j\}}$$
\item {\bf Braid relations:} %
for any connected subdiagram $B\subseteq D$ consisting of two vertices
$\alpha_i,\alpha_j$ such that the order $m_{ij}$ of $s_is_j\in W$ is finite,
the following holds:
\begin{equation}\label{eq:braid relations}
\Ad(\Phi_{(B;\alpha_i,\alpha_j)})(S_i)\cdot S_j\cdots
=
S_j\cdot \Ad(\Phi_{(B;\alpha_i,\alpha_j)})(S_i)\cdots
\end{equation}
where the number of factors on each side is equal to $m_{ij}$.
\end{itemize}
as well as an additional axiom called the {\it generalised pentagon relations},
see \cite[\S 3.17]{TL2}.

The above axioms are designed so that the elements $S_i$ and $\Phi_{(B;
\alpha_i,\alpha_j)}$ define a representation of the Tits braid group $B_W$ 
on any $W$--module, with isomorphic \qca structures yielding equivalent
representations of $B_W$, see \cite[\S 3.14]{TL2}.

\subsection{Deformations of \qca structures \cite[\S 5]{TL2}}

Let now $\R=k\fml$ be the ring of formal power series in a variable $\hbar$.
Let
$$(\{S_i\},\{\Phi_{(B;\alpha_i,\alpha_j)}\})\aand
(\{S_i\pprime\},\{\Phi\pprime_{(B;\alpha_i,\alpha_j)}\})$$
be \qca structures on $\R W=kW\fml$ such that, mod $\hbar$
$$\Phi_{(B;\alpha_i,\alpha_j)}=1=\Phi\pprime_{(B;\alpha_i,\alpha_j)}$$

\noindent
Assume further that $S_i=S_i'$ for any $\alpha_i\in D$, and that the two
structures coincide mod $\hbar^n$ for some $n\geq 1$, that is that
\begin{equation}\label{eq:infinitesimal}
\Phi\pprime_{(B;\alpha_i,\alpha_j)}=
\Phi_{(B;\alpha_i,\alpha_j)}+\hbar^n\varphi_{(B;\alpha_i,\alpha_j)}\mod\hbar^{n+1}
\end{equation}
for any $\alpha_i\neq\alpha_j\in B\subseteq D$, where $\varphi_{(B;\alpha_i,\alpha_j)}
\in{k W_B}^{B\setminus\{\alpha_i,\alpha_j\}}$. 

Then, by \cite[Thm. 5.22]{TL2}, $\varphi=\{\varphi_{(B;\alpha_i,\alpha_j)}\}$
is a 2--cocycle in the Dynkin complex $CD(kW)$ and the two structures
are isomorphic mod $\hbar^{n+1}$ if, and only if, $\varphi$ is a coboundary.

\subsection{Rigidity}

Let $\{S_i\}_{i\in I}$ be elements such that 
$$S_i\in kW_{\alpha_i}\fml\aand S_i=s_i\mod\hbar$$

\begin{theorem}\label{th:rigidity}
Assume that $W$ is finite. Then there exists, up to isomorphism, at most one
\qca structure on $kW\fml$ with local monodromies given by the elements $S
_i$ and associators equal to $1$ mod $\hbar$.
\end{theorem}
\begin{proof}
Let $\varphi\in CD^2(kW)$ be the infinitesimal defined by \eqref{eq:infinitesimal}.
The linearisation of the braid relations \eqref{eq:braid relations} reads, for any
connected subdiagram $B\subseteq D$ with vertex set $\{\alpha_i,\alpha_j\}$
$$\Alt_B(\varphi_{(B;\alpha_i,\alpha_j)})=0$$
where $\Alt_B:kW_B\to{kW_B}^\veps$ is the antisymmetrisation operator.

The image of $\varphi$ in $HC^2(kW)=\bigoplus_{B\subseteq D:|B|=2}{kW_B}
^{\veps}$ via the morphism \eqref{eq:alt quiso} is therefore zero so that $[\varphi]
=0$ in $HD^2(kW)$ by Proposition \ref{pr:Dynkin=Coxeter}.
\end{proof}

\subsection{}

When $k=\IC$, one can endow $kW\fml$ with two \qca structures having
local monodromies 
$$S_i=s_i\cdot \exp(\pi\sqrt{-1}k_{\alpha_i}\hbar s_i)$$
where $q_{\alpha_i}\in\IC$ are a set of complex weights invariant under $W$
\cite[\S 4]{TL2}. The first structure comes from the standard one on the
Iwahori--Hecke algebra $\H W$ obtained by quotienting the group algebra
of the braid group $B_W$ by the quadratic relations
$$(S_i-q_i)(S_i+q_i^{-1})=0$$
where $q_i=\exp(2\pi\sqrt{-1}k_{\alpha_i}\hbar)$, The second one underlies
the monodromy of Cherednik's rational KZ connection \cite{Ch}.

By \cite[\S 4.2.2]{TL2}, these two structures are isomorphic. Theorem \ref
{th:rigidity} strenghtens this result by showing that there are no other such
structures with the above local monodromies.

\section{Top dimensional Dynkin diagram cohomology of Coxeter groups}
\label{sec:top}

\subsection{Sign--coinvariants of $\alg{W}$}
\label{se:signcoinvariants}

Let $W$ be a Coxeter group with system of generators $S=\{s_i\}_{i\in I}$.
Let $V$ be a $W$--module and set
\begin{equation}\label{eq:V_W}
\overline{V}=\sum_i V^{s_i}\aand V_\veps=V/\overline{V}
\end{equation}

\begin{proposition}\label{pr:V_W}\hfill
\begin{enumerate}
\item $\overline{V}$ is invariant under $W$.
\item For any $v\in V$ and $w\in W$,
$$w\medspace v=\veps(w) v\mod\overline{V}$$
where $\veps$ is the sign character of $W$. Thus, $W$ acts on $V_\veps$
by $\veps$.
\item If $V=\bigoplus_j V_j$ is a direct sum of $W$-submodules, then $\overline{V}=
\bigoplus_j \overline{V}_j$. In particular,
$$V_\veps=\bigoplus_j (V_j)_\veps$$
\end{enumerate}
\end{proposition}
\begin{proof}
(i) Since $V^{s_i}=(1+s_i)V$, we have
$$s_j V^{s_i}=
s_j(1+s_i)V\subset (1+s_j)(1+s_i)V+(1+s_i)V\subset V^{s_j}+V^{s_i}$$
(ii) Write $w=s_{i_1}\cdots s_{i_\ell}$. Then, for any $v\in V$
\begin{equation*}
\begin{split}
w v
&=
s_{i_1}\cdots s_{i_\ell} v\\
&=
-s_{i_2}\cdots s_{i_\ell}v + (1+s_{i_1})s_{i_2}\cdots s_{i_\ell}v\\
&=
(-1)^\ell v+\sum_{j=1}^\ell (-1)^{j-1}(1+s_{i_j})s_{i_{j+1}}\cdots s_{i_\ell}v
\end{split}
\end{equation*}
(iii) is clear.
\end{proof}

\remark Proposition \ref{pr:V_W} also follows from the fact that $V^{s_i}=
(1+s_i)V$ so that $V_\veps\simeq \veps\otimes_{kW}V$.
\smallskip

Let now $V=\alg{W}$ endowed with the conjugation action of $W$ and let
$\CC$ be the set of conjugacy classes of $W$. For any $c\in\CC$, choose
$w_c\in c$ and let $A_{w_c}$ be the image of $w_c\in\alg{W}$ in $\alg{W}
_\veps$. Let $C_W(w_c)$ be the centraliser of $w_c$ in $W$.

\begin{proposition}\label{pr:sign}\hfill
\begin{enumerate}
\item For any Coxeter group $W$, we have $A_{w_c}\not=0$ when 
$\veps(C_W(w_c))=1$ and
$$\alg{W}_\veps=
\bigoplus_{\substack{c\in\CC:\\ \veps(C_W(w_c))=1}}kA_{w_c}$$
\item If $W$ is finite, then given $c$ such that
$\veps(C_W(w_c))=1$, the element
$A^{w_c}=\sum_{w'\in W/C_W(w_c)}\veps(w') w'w_cw^{\prime -1}$ is non-zero
and
$$\alg{W}^\veps=
\bigoplus_{\substack{c\in\CC:\\ \veps(C_W(w_c))=1}}kA^{w_c}.$$
\end{enumerate}
\end{proposition}
\begin{proof}
(i) Since $\alg{W}=\bigoplus_{c\in\CC}\CF_c$, where $\CF_c=\vec{c}$ is the
subspace spanned by elements of $c$, Proposition \ref{pr:V_W} yields $
\alg{W}^\veps_W=\bigoplus_{c\in\CC}(\CF_c)_\veps$. Since $W$ acts
transitively on $c$, it follows from (ii) of Proposition \ref{pr:V_W} that $
(\CF_c)_\veps$ is spanned by $A_{w_c}$ and therefore at most one--dimensional.
If $w$ centralises $w_c$, then, by (ii) of Proposition \ref{pr:V_W},
$$A_{w_c}=A_{ww_cw^{-1}}=w A_{w_c}=\veps(w) A_{w_c}$$
so that $A_{w_c}$ is zero if the sign character is not trivial on the centraliser
of $w_c$. Conversely, if $\veps(C_W(w_c))=1$, the assignement ${ww_c}
\rightarrow\veps(w)$ consistently defines a non--zero linear form on $\CF
_c$ which descends to $(\CF_c)_\veps$ so that $A_{w_c}\neq 0$. (ii)
readily follows from (i) and the fact that if $W$ is finite, $V_\veps\simeq
V^\veps$ for any $W$--module $V$.
\end{proof}

\medskip
Consider $c\in\CC$ such that $\veps(C_W(w_c))=1$. Then, $A_{w_c}$
depends only on $c$ and we put $A_c=A_{w_c}$. Similarly, when $W$ is finite
we put $A^c=A^{w_c}$.




\subsection{Top cohomology}

Assume now that $W$ is irreducible of rank $n$ and let $D$ be its Coxeter
graph. Let $\CP$ be the collection of proper, maximal connected subdiagrams
$B$ of $D$ and $\alg{W}_\CP\subset \alg{W}$ the span of $\alg{W_B}$ as
$B$ varies in $\CP$. For any $c\in\CC$, let $\overline{A}_c$ be the class
of $A_c$ in $\left.\alg{W}_\veps\right/ \left(\alg{W}_\CP\cap \overline{kW}
\right)$.

\begin{proposition}\label{pr:top Dynkin}
We have a decomposition in one-dimensional subspaces
$$HD^n(\alg{W})=
\bigoplus_{\substack
{c\in\CC:\\ \veps(C_W(w_c))=1,\\ c\cap W_B=\emptyset,\;\forall B\in\CP}}
k\overline{A}_{w_c}$$
\end{proposition}
\begin{proof}
The top degree Dynkin differential $d_D^n$ is zero. Moreover,
$$d_D^{n-1}a_{(D;D)}=
\sum_{i=1}^n (-1)^{i-1}
(a_{(D;D\setminus\alpha_i)}-
a_{(\con{D}{\alpha_i}{D\setminus\alpha_i};D\setminus\alpha_i)})$$
Since $a_{(D;D\setminus\alpha_i)}\in \alg{W}^{s_i}$ and $a_{(\con{D}
{\alpha_i}{D\setminus\alpha_i};D\setminus\alpha_i)}\in\alg{W_{\con{D}
{\alpha_i}{D\setminus\alpha_i}}}$,
$$\Im d_D^{n-1}=\sum_{i=1}^n \alg{W}^{s_i}+\sum_{B\in\CP}\alg{W_B}$$
The result now follows from Proposition \ref{pr:sign}.
\end{proof}

\section{Finite Coxeter groups}
\label{sec:finite}

\subsection{Finite Coxeter groups of rank 2}\label{ss:rank 2}

Let $W=\sfI_2(m)$, $m\geq 3$, be the Coxeter group with generators
$s,t$ and relations $s^2=1=t^2$ and $(st)^m=1$. For $p=0,\ldots,m-1
$, let $c^p$ be the conjugacy class of $(st)^p$. 

\begin{proposition}\label{pr:rank 2 sign}
$$\alg{W}^\veps=\bigoplus_{p=1}^{\lfloor\frac{m-1}{2}\rfloor} k A^{c^p}$$
\end{proposition}
\begin{proof}
The only conjugacy classes involved in the decomposition of proposition
\ref{pr:sign} are those of words in $s,t$ of even length and therefore
those of the powers of $st$ and $ts$. Since $ts=s(st)s$ we need
only consider the classes $c^p$, $p=1,\ldots,m-1$. Moreover, since
$(st)^{m-p}=(ts)^p=s(st)^ps$, $c^p=c^{m-p}$ and we may restrict
our attention to $1\leq p\leq m/2$. Finally, since for $m$ even $(st)^
{m/2}$ is central in $W$, the only possible relevant values of $p$
are $1,\ldots,\lfloor\frac{m-1}{2}\rfloor$. The proposition follows from the
fact that for $p=1,\ldots,\lfloor\frac{m-1}{2}\rfloor$, the centraliser of
$(st)^p$ in $W$ is generated by $st$.
\end{proof}

Since the differential $d^\#$ \eqref{eq:d flat} is zero for $W$ of rank
2, Proposition \ref{pr:rank 2 sign} implies the following

\begin{theorem}
For any $m\geq 3$,
$$\dim HD^p(\alg{\sfI_2(m)})=
\left\{\begin{array}{cl}
0 &\text{if $p=0$}\\
0 &\text{if $p=1$}\\
\lfloor\frac{m-1}{2}\rfloor &\text{if $p=2$}\\
0 &\text{if $p\geq 3$}
\end{array}\right.$$
\end{theorem}

\subsection{Type $\sfA_n$}\label{se:An}

Let $W=\GS_{n+1}$ be the Weyl group of type $\sfA_n$, $n\geq 1$.
The conjugacy classes in $W$ are parametrised by partitions $\lambda
=(\lambda_1,\ldots,\lambda_k)$ of $n+1$, with $c^\lambda$ the class
of the product of cycles
$$\tau^\lambda=
(1\ths 2\ths\cdots\ths\lambda_1)
(\lambda_1+1\ths\lambda_1+2\ths\cdots\ths\lambda_1+\lambda_2)
\cdots
(\sum_{i =1}^{k-1}\lambda_i+1\ths\cdots\ths\sum_{i=1}^k\lambda_i)$$

\subsubsection{} 

For any $m\in\IN^*$, let
\begin{equation}\label{eq:odd}
\CO_m=
\{\lambda\vdash m|\thinspace
\lambda_i\in 2\IN+1,\forall i
\medspace\text{and}\medspace
\lambda_i\neq\lambda_j,\thinspace\forall i\neq j\}
\end{equation}
be the set of partitions of $m$ consisting of odd, distinct parts.

\begin{proposition}\label{pr:A sign}
$$\alg{\GS_{n+1}}^\veps=
\bigoplus_{\lambda\in\CO_{n+1}}kA^{c^\lambda}$$
\end{proposition}
\begin{proof}
Since $\veps(m\ths m+1\ths\cdots\ths m+p-1)=(-1)^{p-1}$, the only
conjugacy classes involved in the decomposition of Proposition \ref{pr:sign}
are those such that each $\lambda_i$ is odd. Moreover, since the product 
$$(m\ths m+1\ths\cdots\ths m+p-1)(m'\ths m'+1\ths\cdots\ths m'+p-1)$$
of two disjoint cycles of equal length is centralised by $\pi=(m\ths m')
\cdots(m+p-1\ths m'+p-1)$ and $\veps(\pi)=(-1)^p$, the $\lambda_i$'s
must all be distinct. When this last condition is fulfilled, any element
centralising $\tau^\lambda$ is of the form
$$\zeta=
(1\ths 2\ths\cdots\ths\lambda_1)^{m_1}
(\lambda_1+1\ths\lambda_1+2\ths\cdots\ths\lambda_1+\lambda_2)^{m_2}
\cdots
(\sum_{i =1}^{k-1}\lambda_i+1\ths\cdots\ths\sum_{i=1}^k\lambda_i)^{m_k}$$
for some $0\leq m_i\leq\lambda_i-1$. Since $\veps(\zeta)=\prod_i
((-1)^{\lambda_i-1})^{m_i}$, the partitions arising in the decomposition
of Proposition \ref{pr:sign} are exactly those in $\CO_{n+1}$.
\end{proof}

\subsubsection{} 

Identify the connected subdiagrams of the Coxeter graph $D$ of $W$ with
the subintervals of $[1,n]$ having integral endpoints so that $W_{[i,j]}\simeq
\GS_{j-i+2}$. For any $2\leq p\leq n$ and $1\leq i\leq n-p+1$, let
\begin{equation}\label{eq:Aclambda}
A^{c^\lambda}_{[i,i+p-1]}
=\Alt_{[i,i+p-1]}\left(\tau_{[i,i+p-1]}^\lambda\right)
\in kW_{[i,i+p-1]}^\veps
\end{equation}
be the generator corresponding to $\lambda\in\CO_{p+1}$. We shall need
the following

\begin{lemma}\label{le:A in A}
\begin{align*}
\Alt_{[i,i+p]}(A_{[i,i+p-1]}^{c^\lambda})
&=\left\{\begin{array}{cl}
A_{[i,i+p]}^{c^{\lambda\cupone}}&\text{if $1\notin\lambda$}\\
0&\text{otherwise}
\end{array}\right.
\intertext{and}
\Alt_{[i-1,i+p-1]}(A_{[i,i+p-1]}^{c^\lambda})
&=\left\{\begin{array}{cl}
(-1)^{p+1} A_{[i-1,i+p-1]}^{c^{\lambda\cupone}}&\text{if $1\notin\lambda$}\\
0&\text{otherwise}
\end{array}\right.
\end{align*}
provided $i\leq n-p$ and $i\geq 2$ respectively.
\end{lemma}
\begin{proof}
The first identity follows from the fact that under the inclusion $\GS_{[i,i+p-1]}
\subset\GS_{[i,i+p]}$, $\tau_{[i,i+p-1]}^\lambda$ is mapped to $\tau_{[i,i+p]}^
{\lambda\cupone}$. The second one follows from the first and the fact that,
in $\GS_{[i-1,i+p-1]}$,
$$\tau_{[i,i+p-1]}^\lambda=
\Ad((i-1\ths i\ths\cdots\ths i+p))\tau_{[i-1,i+p-2]}^\lambda$$
and $\veps(i-1\ths i\ths\cdots\ths i+p)=(-1)^{p+1}$.
\end{proof}

\subsubsection{} 

Label the nodes of $D$ as in \cite[Planche I]{Bo} and order them as $\alpha_1
<\cdots<\alpha_n$. For any $p=1,\ldots,n$, let
$$d^{\#}_p:
\bigoplus_{i=1}^{n-p+1}\alg{W_{[i,i+p-1]}}^\veps\longrightarrow
\bigoplus_{j=1}^{n-p}\alg{W_{[j,j+p]}}^\veps$$
be the differential of $HC$, where the \rhs is understood to
be $0$ if $p=n$. Set
\begin{equation}\label{eq:O*}
\CO_{p+1}^*=\{\lambda\in\CO_{p+1}|\ths 1\notin\lambda\}
\end{equation}

\begin{proposition}\label{pr:df for A}
For any $p=2,\ldots,n$,
\begin{align}
\Im d^\#_{p-1}
&=
\bigoplus_{i=1}^{n-p+1}
\bigoplus_{\substack{\lambda\in\CO_{p+1}:\\1\in\lambda}}
k A^{c^\lambda}_{[i,i+p-1]}
\label{eq:Im}\\
\Ker d^\#_p
&=
\Im d^\#_{p-1}\oplus
\bigoplus_{\lambda\in\CO_{p+1}^*}
k B^\lambda
\label{eq:Ker}
\end{align}
where $B^\lambda=\sum_{i=1}^{n-p+1}A^{c^\lambda}_{[i,i+p-1]}$.
\end{proposition}
\begin{proof}
Since
$$([i-1,i+p-2];[i,i+p-2])=1\quad\text{and}\quad([i,i+p-1];[i,i+p-2])=p$$
Lemma \ref{le:A in A} yields, for any $\overlineambda\in\CO_p$ and $i=1,\ldots,n-p+2$,
\begin{multline}\label{eq:df A}
d^\#_{p-1} A_{[i,i+p-2]}^{c^\overlineambda}\\
=
\delta_{1\notin\overlineambda}\cdot(-1)^p\left(
-\delta_{i\geq 2}\cdot A_{[i-1,i+p-2]}^{c^{\overlineambda\cupone}}
+\delta_{i\leq n-p+1}\cdot A_{[i,i+p-1]}^{c^{\overlineambda\cupone}}
\right)
\end{multline}
It follows that, for any $\lambda\in\CO_p^*$,
\begin{equation}\label{eq:coboundary}
d^\#_{p-1}(\sum_{j=1}^i A_{[j,j+p-2]}^{c^\overlineambda})=
(-1)^p A_{[i,i+p-1]}^{c^{\overlineambda\cupone}}
\end{equation}
which yields \eqref{eq:Im}. It also follows from \eqref{eq:df A} that, for
$\lambda\in\CO_{p+1}^*$, $B^\lambda$ is the unique linear combination
of $A_{[i,i+p-1]}^{c^\lambda}$ such that $d^\#_p B^\lambda=0$, which
yields \eqref{eq:Ker}.
\end{proof}

\subsubsection{} 

\begin{theorem}
For any $0\leq p\leq n$,
\begin{equation}\label{eq:HDA}
\dim HD^p(k\GS_{n+1})_{c^\lambda}=
\begin{cases}
1 & \text{ if } \lambda=\lambda'\cup \{1^{n-p}\} \text{ and }
\lambda'\in \CO_{p+1}^*\\
0 & \text{ otherwise}
\end{cases}
\end{equation}
\end{theorem}
\begin{proof}
\eqref{eq:HDA} holds for $p\geq 2$ by proposition \ref{pr:df for A}
and therefore for $p\geq 0$ since $\dim HD^i(k\GS_{n+1})=0=|\CO
_{i+1}^*|$ for $i=0,1$.
\end{proof}

\subsubsection{Generating function}

\begin{theorem}\label{th:A}
Set
$$\chi^{\sfA}(q,t)=
\sum_{n\geq 1,p\geq 0}
q^n t^p\dim HD^p(\alg{\GS_{n+1}})$$
Then,
$$\chi^{\sfA}(q,t)=
\frac{1}{1-q}\frac{\prod_{d\geq 1}(1+(qt)^{2d+1})-1}{qt}$$
\end{theorem}
\begin{proof}
Since $HD^p(k\GS_{n+1})=0$ if $p=0$ or $p\geq n+1$, we have
\begin{equation*}
\begin{split}
\chi^{\sfA}(q,t)
&=
\sum_{n\geq p\geq 1} q^n t^p\dim HD^p(\alg{\GS_{n+1}})\\
&=
\sum_{p\geq 1}\frac{t^p q^p}{1-q}|\CO_{p+1}^*|\\
&=
\frac{1}{1-q}\frac{\prod_{d\geq 1}(1+(qt)^{2d+1})-1}{qt}
\end{split}
\end{equation*}
where the last identity follows from the fact that
\begin{equation}\label{eq:phi*}
\sum_{m\geq 1}z^m|\CO_m^*|=
\sum_{m\geq 2}z^m|\CO_m^*|=
(1+z^3)(1+z^5)(1+z^7)\cdots - 1
\end{equation}
\end{proof}

\remark Up to a multiplication by $qt$, the generating function
$\chi^{\sfA}$ is the product a Fermionic partition function by
that of a one--dimensional harmonic oscillator. It would be
interesting to know whether the direct sum
$$\bigoplus_{n\geq 1,p\geq 0}HD^p(k\GS_{n+1})$$
possesses a natural action of an infinite--dimensional Clifford
algebra similar in spirit to that on the cohomology of the Hilbert
schemes of points on a surface \cite{Gr,Na}.

\subsection{Type $\sfB_n$}\label{se:Bn}

Let now $W=\GS_n\ltimes\IZ_2^n$ be the Weyl group of type $\sfB
_n$, $n\geq 2$, and denote the generators of $\IZ_2^n$ by $\veps
_i$, $i=1,\ldots,n$. The conjugacy classes in $W$ are parametrised
by ordered pairs of partitions $(\lambda,\mu)$ such that $|\lambda|
+|\mu|=n$, where $|\lambda|=\sum_i\lambda_i$ \cite[Prop. 3.4.7]{GP}.
The class $c^{(\lambda,\mu)}$ is that of the product $\tau^\lambda
\wt{\tau}^\mu$, where
\begin{align}
\tau^\lambda
&=
(1\ths 2\ths\cdots\ths\lambda_1)
(\lambda_1+1\ths\lambda_1+2\ths\cdots\ths\lambda_1+\lambda_2)
\cdots
(|\lambda|-\lambda_k+1\ths\cdots\ths|\lambda|)
\label{eq:tau lambda}\\
\intertext{and}
\wt{\tau}^\mu
&=
(|\lambda|+1\ths\cdots\ths|\lambda|+\mu_1)\veps_{|\lambda|+\mu_1}
\cdots
(n-\mu_\ell+1\ths\cdots\ths n)
\veps_n
\label{eq:tau mu}
\end{align}

\begin{proposition}\label{pr:B sign}
$$\vec{\GS_n\ltimes\IZ_2^n}^\veps=
\left\{\begin{array}{cl}
\displaystyle{\bigoplus_{\nu\vdash n/2}
kA^{c^{(\emptyset,2\nu)}}}
&\text{if $n$ is even}\\[1.2ex]
0&\text{if $n$ is odd}
\end{array}\right.$$
\end{proposition}
\begin{proof}
A necessary condition for a conjugacy class $c^{(\lambda,\mu)}$ to
contribute to the decomposition of Proposition \ref{pr:sign} is that
the $\lambda_i$ be odd and the $\mu_j$ even since $\veps(\veps_i)=-1$.
Since $(m\ths m+1\ths\cdots\ths m+p-1)$ is centralised by $\veps_m
\cdots\veps_{m+p-1}$ and $\veps(\veps_m\cdots\veps_{m+p-1})=(-1)
^p$, $\lambda$ must in fact be the empty partition. In particular,
$n$ must be even\footnote{this also follows from the fact that the
central element $\zeta=\veps_1\cdots\veps_n$ is of sign $(-1)^n$ so
that $$\alg{W}^\veps\subseteq\{f\in \alg{W}|\medspace \zeta f\zeta=(-1)^n f\}$$
is zero if $n$ is odd.}. There remains to show that if $\mu\vdash n$
only contains even parts, the sign character $\veps$ is trivial on the
centraliser of $\wt{\tau}^\mu$. This is readily reduced to the case
when $\mu$ only has one part which follows in turn from the following
result.
\end{proof}

\begin{lemma}\label{le:signed cycle}
The centraliser of $(1\ths\cdots\ths p)\veps_p$ in $\GS_p\ltimes
\IZ_2^p$ is the group $\IZ_{2p}$ generated by $(1\ths\cdots\ths
p)\veps_p$.
\end{lemma}
\begin{proof}
If $w\in\GS_p\ltimes\IZ_2^p$ centralises $(1\ths\cdots\ths
p)\veps_p$, its projection in $\GS_p$ centralises $(1\ths\cdots
\ths p)$ and is therefore equal to $(1\ths\cdots\ths p)^s$ for 
some $0\leq s\leq p-1$. Thus $((1\ths\cdots\ths p)\veps_p)^{-s}
w$ centralises $(1\ths\cdots\ths p)\veps_p$ and lies in $\IZ_2^p$
from which it readily follows that it is either equal to 1 or
to $\veps_1\cdots\veps_p=((1\ths\cdots\ths p)\veps_p)^p$.
\end{proof}

\begin{proposition}\label{pr:HDB}
For any $0\leq p\leq n$, we have
$$\dim HD^p(kW_{B_n})_{c^{\lambda,\mu}}=
\begin{cases}
1 & \text{ if } \left(\mu=0 \text{ and }
\lambda=\lambda'\cup \{1^{n-p-1}\} \text{ with }
\lambda'\in \CO_{p+1}^*\right) \\
 & \text{ or } \left(\lambda=(1^{n-p}) \text{ and }
 \mu=2\nu \text{ with }\nu\vdash p/2
\right)\\
0 & \text{ otherwise.}
\end{cases}
$$
In particular,
$$\dim HD^p(W_{B_n})=
\dim HD^p(\alg{\GS_n})+
\delta_{p\in 2\IN^*}\cdot P(p/2)$$
where $P$ is the partition function.
\end{proposition}
\begin{proof}
Identify the connected subdiagrams of the Coxeter graph $D$
of $W$ with the subintervals of $[1,n]$ having integral endpoints so
that
$$W_{[i,j]}\simeq
\left\{\begin{array}{ll}
\GS_{j-i+2}                    &\text{if $j\leq n-1$}\\
\GS_{n-i+1}\ltimes\IZ_2^{n-i+1}&\text{if $j=n$}
\end{array}\right.$$
For $1\leq p\leq n-1$ and $1\leq i\leq n-p$, let
$$A_{[i,i+p-1]}^{c^\lambda}\in kW_{[i,i+p-1]}^\veps$$
be the generator corresponding to $\lambda\in\CO_{p+1}$, as in \eqref
{eq:Aclambda}. By Proposition \ref{pr:B sign}, $\Alt_{[i,n]}(A_{[i,n-1]}^
{c^\lambda})=0$. Since, in addition $\alg{W_{[n-p+1,n]}}^\veps$ is zero
whenever $p$ is odd and of dimension $|\{\nu\vdash p/2\}|=P(p/2)$
when $p$ is even, the complex $\left(\bigoplus_{B\subseteq D} \alg{W_B}
^\veps;d^\#\right)$ of \S \ref{se:E2} decomposes as the direct sum of the
corresponding complex for $\GS_n$ and a complex concentrated in
positive, even degrees with chain groups of dimension $\delta_{p\in
2\IN^*}\cdot P(p/2)$.
\end{proof}

\begin{theorem}\label{th:B}
Set $\chi^{\sfB}(q,t)=\sum_{n\geq 2,p\geq 0}q^n t^p\dim HD^p
(kW_{B_n})$. Then,
$$\chi^{\sfB}(q,t)=
\frac{\prod_{d\geq 1}(1+(qt)^{2d+1})-1}{(1-q)t}+
\frac{\prod_{d\geq 1}(1-(qt)^{2d})^{-1}-1}{1-q}$$
\end{theorem}
\begin{proof}
This follows from Proposition \ref{pr:HDB}, Theorem \ref{th:A}
and the fact that
$$\sum_{m\geq 0}z^m P(m)=\prod_{d\geq 1}(1-z^d)^{-1}$$
\end{proof}

\subsection{Type $\sfD_n$}\label{se:Dn}

Let $\IZ\tp^n\subset\IZ_2^n$ be the kernel of the sign character and
$W=\GS_n\ltimes\IZ\tp^n$ the Weyl group of type $\sfD_n$, $n\geq 3$.
The conjugacy classes in $W$ fall into two types \cite[Prop. 3.4.12]{GP}:
\begin{enumerate}
\item[type I.] These are labelled by ordered pairs $(\lambda,\mu)$ of
partitions such that $|\lambda|+|\mu|=n$ and the number of parts
$[\mu]$ of $\mu$ is even. The corresponding class $c^{(\lambda,
\mu)}$ is that of the product $\tau^\lambda\wt{\tau}^\mu$, where
$\tau^\lambda$, $\wt{\tau}^\mu$ are given by
\eqref{eq:tau lambda}--\eqref{eq:tau mu}.
\item[type II.] These are labelled by partitions $\lambda$ of $n$
all of whose parts are even, with $c^{\lambda,\II}$ the class of
$$\tau^{\lambda,\II}=
(1\ths 2\ths\cdots\ths\lambda_1)
\cdots
(n-\lambda_{k-1}-\lambda_k+1\ths\cdots\ths n-\lambda_k)
(n-\lambda_k+1\ths\cdots\ths n)\veps_{n-1}\veps_n$$
\end{enumerate}

\subsubsection{} 

\begin{proposition}\label{pr:D sign}
\begin{equation}\label{eq:D sign}
(k\GS_n\ltimes\IZ\tp^n)^\veps=
\bigoplus_{\substack
{0\leq m\leq n,\\[.2ex]\lambda\in\CO_m,\\\mu\in\CO_{n-m}:\\ [\mu]\in 2\IN}}
kA^{c^{(\lambda,\mu)}}
\oplus\left\{\begin{array}{cl}
\displaystyle{\bigoplus_{\substack{\nu\vdash n/2:\\ [\nu]\in 2\IN}}
kA^{c^{(\emptyset,2\nu)}}}
&\text{if $n$ is even}\\
0&\text{if $n$ is odd}
\end{array}\right.
\end{equation}
\end{proposition}
\begin{proof}
Conjugacy classes of type II do not contribute to the decomposition
of Proposition \ref{pr:sign} since $\tau^{\lambda,\II}$ is centralised
by
$$w=(n-\lambda_k+1\ths\cdots\ths n)\veps_{n-1}\veps_n$$
and $\veps(w)=(-1)^{\lambda_k-1}=-1$. Consider now a class
$c^{(\lambda,\mu)}$ of type I which contributes to the decomposition
of Proposition \ref{pr:sign}. As in Proposition \ref{pr:A sign}, the
$\lambda_i$ must be odd and distinct if $\lambda$ is nonempty.
Moreover, since for any $1\leq i<j\leq\ell$
\begin{equation*}
\begin{split}
w_{i,j}
=&
(|\lambda|+\sum_{a=1}^{i-1}\mu_a+1\ths\cdots\ths |\lambda|+\sum_{a=1}^i\mu_a)
\veps_{|\lambda|+\sum_{a=1}^i\mu_a}\\
\cdot&
(|\lambda|+\sum_{a=1}^{j-1}\mu_a+1\ths\cdots\ths |\lambda|+\sum_{a=1}^j\mu_a)
\veps_{|\lambda|+\sum_{a=1}^j\mu_a}\in W
\end{split}\end{equation*}
centralises $\tau^\lambda\wt{\tau}^\mu$ and $\veps(w_{i,j})=(-1)^
{\mu_i+\mu_j}$, all $\mu_i$ must be of the same parity. Finally,
since for $\lambda_1$ odd,
$$w=(1\ths 2\ths\cdots\ths\lambda_1)\veps_1\cdots\veps_{\lambda_1}
\cdot
(|\lambda|+1\ths\cdots\ths|\lambda|+\mu_1)\veps_{|\lambda|+\mu_1}$$
lies in $W$, centralises $\tau^\lambda\wt{\tau}^\mu$ and $\veps(w)=(-1)
^{\mu_1-1}$, all $\mu_i$ must be odd, and therefore distinct, if $\lambda$
is nonempty.

There remains to show that if $c$ is a conjugacy class appearing on
the \rhs of \eqref{eq:D sign}, then $\veps$ is trivial on the centraliser
of any element of $c$. If $c=c^{(\emptyset,\mu)}$ with all $\mu_i$
even, the centraliser of $\wt{\tau}^\mu$ in $\IZ_2^n\rtimes\GS_n$
lies in the kernel of $\veps$ by Proposition \ref{pr:B sign}. {\it A fortiori},
this is true in $\IZ\tp^n\rtimes\GS_n$. If, on the other hand, $c=c^
{(\lambda,\mu)}$, where $\lambda$ and $\mu$ are either empty
or consist of odd, distinct parts, it follows from Lemma \ref{le:centraliser}
below that the component $\sigma\in\GS_n$ of any $w\in W$ centralising
$\tau^\lambda\wt{\tau}^\mu$ is of the form
\begin{equation*}
\begin{split}
&
(1\ths 2\ths\cdots\ths\lambda_1)^{m_1}
\cdots
(|\lambda|-\lambda_k+1\ths\cdots\ths|\lambda|)^{m_k}\\
\cdot&
(|\lambda|+1\ths\cdots\ths|\lambda|+\mu_1)^{m'_1}
\cdots
(n-\mu_\ell+1\ths\cdots\ths n)^{m'_\ell}
\end{split}
\end{equation*}
so that $\veps(w)=\veps(\sigma)=1$.
\end{proof}

\begin{lemma}\label{le:centraliser}
The centraliser of $(1\cdots p)(p+1\cdots 2p)\veps_{2p}$ in in $\GS_{2p}
\ltimes\IZ_2^{2p}$ is the product of the centralisers of $(1\cdots p)$ and
$(p+1\cdots 2p)\veps_{2p}$ in $\GS_{\{1,\ldots,p\}}\ltimes\IZ_2^p$ and
$\GS_{\{p+1,\ldots,2p\}}\ltimes\IZ_2^p$ respectively.
\end{lemma}

\subsubsection{} 

Label now the nodes of the Dynkin diagram $D$ of $W$ as in \cite[Planche IV]
{Bo}, so that $\alpha_{n-2}$ is the trivalent node of $D$ if $n\geq 4$, and order
them as $\alpha_1<\cdots<\alpha_{n-2}<\alpha_{n-1}<\alpha_n$. For any
$1\leq i\leq j\leq n-2$, let $[i,j]\subset D$ be the connected subdiagram
with vertices $\alpha_i,\ldots,\alpha_j$. For any $i=1,\ldots,n-2$, let $B_i
^\pm\subset D$ be the connected subdiagrams with vertices $\alpha_i,
\ldots,\alpha_{n-2}$ and $\alpha_{n-1}$ (resp.  $\alpha_n$) and $B_i^Y
\subseteq D$ the subdiagram with vertices $\alpha_i,\ldots,\alpha_n$.
Thus,
$$W_{[i,j]}\simeq\GS_{j-i+2},\quad
W_{B_i^\pm}\simeq\GS_{n-i+1}\quad\text{and}\quad
W_{B_i^Y}\simeq\GS_{n-i+1}\ltimes\IZ\tp^{n-i+1}$$

Let $\sigma\in\Aut(W)$ be the involution induced by fixing the nodes
$\alpha_1,\ldots,\alpha_{n-2}$ and permuting $\alpha_{n-1}$ and $
\alpha_n$ so that $\sigma(W_{B_i^+})=W_{B_i^-}$. For any $\lambda\in
\CO_{n-i+1}$, set $A_{B_i^-}^{c^\lambda}=\sigma A_{B_i^+}^{c^\lambda}$.


\begin{proposition}\label{pr:Im D}
The following holds for $3\leq p\leq n$,
\begin{equation}
\begin{split}
\Im d^\#_{p-1}=
& 
\bigoplus_{\substack{\lambda\in\CO_{p+1}:\\1\in\lambda}}
\bigoplus_{i=1}^{n-p-1}
k A_{[i,i+p-1]}^{c^\lambda}\\
\oplus& 
\bigoplus_{\substack{\lambda\in\CO_{p+1}:\\1\in\lambda}}
V^\lambda\\
\oplus& 
\phantom{\negmedspace\negmedspace}
\bigoplus_{\substack{1\leq m\leq p,\\\overlineambda\in\CO_m,\\\overline{\mu}\in\CO_{p-m}:
\\1\in\overlineambda,[\overline{\mu}]\in 2\IN}}
k A_{B_{n-p+1}^Y}^{c^{(\overlineambda,\overline{\mu})}}
\end{split}
\label{eq:Im D}
\end{equation}
where the first summand only arises if $p\leq n-2$ and
\begin{equation}\label{eq:V lambda}
V^\lambda=\left\{\begin{array}{cl}
\displaystyle{
\langle A_{B^+_{n-p}}^{c^\lambda}+A_{B^-_{n-p}}^{c^\lambda},
A_{B^+_{n-p}}^{c^\lambda}-A_{B_{n-p+1}^Y}^{c^{(\lambda\setminus\{1\},\emptyset)}},
A_{B^-_{n-p}}^{c^\lambda}+A_{B_{n-p+1}^Y}^{c^{(\lambda\setminus\{1\},\emptyset)}}\rangle}
&\text{if $p\leq n-1$}\\[2.5ex]
\displaystyle{k A_{B_1^Y}^{c^{(\lambda\setminus\{1\},\emptyset)}}}
&\text{if $p=n$}
\end{array}\right.
\end{equation}
\end{proposition}
\begin{proof}
By \eqref{eq:coboundary}, the image of the restriction of
$d^\#_{p-1}$ to $$\bigoplus_{j=1}^{n-p-1}\alg{W_{[j,j+p-2]}}^\veps$$
is the span of the generators $A_{[i,i+p-1]}^{c^\lambda}$, $i=1,
\ldots,n-p-1$, as $\lambda$ runs through the elements of $\CO_
{p+1}$ containing 1. This accounts for the first summand in \eqref
{eq:Im D}. Further, since
$$([n-p-1,n-2];[n-p,n-2])=1
\quad\text{and}\quad
(B_{n-p}^\pm;[n-p,n-2])=p$$
Lemma \ref{le:A in A} yields, for any $\overlineambda\in\CO_p$,
$$d^\#_{p-1}(A_{[n-p,n-2]}^{c^\overlineambda})=
(-1)^p\cdot\delta_{1\notin\overlineambda}\cdot\left(
-A_{[n-p-1,n-2]}^{c^{\overlineambda\cupone}}+
A_{B_{n-p}^+}^{c^{\overlineambda\cupone}}+
A_{B_{n-p}^-}^{c^{\overlineambda\cupone}}\right)$$
Thus, if $p\leq n-1$ and $\lambda\in\CO_{p+1}$ contains $1$,
\begin{equation}\label{eq:B+ B-}
A_{B_{n-p}^+}^{c^\lambda}+A_{B_{n-p}^-}^{c^\lambda}
\in\Im d^\#_{p-1}
\end{equation}
To proceed, we need the following

\begin{lemma}\label{le:A in D}
For any $\overlineambda\in\CO_p$,
$$\Alt_{B_{n-p+1}^Y}(A_{B_{n-p+1}^\pm}^{c^\overlineambda})=
A_{B_{n-p+1}^Y}^{c^{(\overlineambda,\emptyset)}}$$
\end{lemma}
\begin{proof} The "+" identity follows from the fact that, under the inclusion
$W_{B_{n-p+1}^+}\subset W_{B_{n-p+1}^Y}$, $\tau^\overlineambda_
{B_{n-p+1}^+}$ is mapped to $\tau^\overlineambda_{B_{n-p+1}^Y}$.
The "-" one follows by applying the automorphism $\sigma$ and
noticing that $A_{B_{n-p+1}^Y}^{c^{(\overlineambda,\emptyset)}}$ is
fixed by $\sigma$. Indeed if $1\in\overlineambda$, $\tau^\overlineambda_
{B^Y_{n-p+1}}$ lies in $W_{[n-p+1,n-2]}$ and is therefore fixed
by $\sigma$. If on the other hand $1\notin\overlineambda$, the cycle
$(n-\overlineambda_k+1\ths\cdots\ths n)$ is the product $s_{n-\overlineambda
_k+1}\cdots s_{n-1}$, so that
\begin{equation*}
\begin{split}
\sigma(n-\overlineambda_k+1\ths\cdots\ths n)
&=
s_{n-\overlineambda_k+1}\cdots s_{n-2}s_n\\
&=
(n-\overlineambda_k+1\ths\cdots\ths n)\veps_{n-1}\veps_n\\
&=
\Ad(\veps_{n-\overlineambda_k+1}\cdots\veps_{n-1})
(n-\overlineambda_k+1\ths\cdots\ths n)
\end{split}
\end{equation*}
whence
$$\sigma\medspace A_{B_{n-p+1}^Y}^{c^{(\overlineambda,\emptyset)}}=
\Alt_{B_{n-p+1}^Y}\left(\Ad(\veps_{n-\overlineambda_k+1}\cdots\veps_{n-1})
\tau^\overlineambda_{B^Y_{n-p+1}}\right)=
A_{B_{n-p+1}^Y}^{c^{(\overlineambda,\emptyset)}}$$
since $\veps(\veps_{n-\overlineambda_k+1}\cdots\veps_{n-1})=
(-1)^{\overlineambda_k-1}$ and $\overlineambda_k$ is odd.
\end{proof}

Since
$$(B_{n-p}^\pm;B_{n-p+1}^\pm)=1,\thickspace
(B_{n-p+1}^Y;B_{n-p+1}^+)=p
\thickspace\text{and}\thickspace
(B_{n-p+1}^Y;B_{n-p+1}^-)=p-1$$
lemmas \ref{le:A in A} and \ref{le:A in D} imply that, for any $\overlineambda
\in\CO_p$,
\begin{equation}\label{eq:d on B+-}
d^\#_{p-1}(A_{B_{n-p+1}^\pm}^{c^\overlineambda})
=(-1)^p\left(
-\delta_{1\notin\overlineambda}\cdot A_{B_{n-p}^\pm}^{c^{\overlineambda\cupone}}
\pm A_{B_{n-p+1}^Y}^{c^{(\overlineambda,\emptyset)}}\right)
\end{equation}
Choosing $\overlineambda\in\CO_p$ such that $1\notin\overlineambda$ in \eqref
{eq:d on B+-} and using \eqref{eq:B+ B-} accounts for the second
summand in \eqref{eq:Im D}. To conclude, we need the following
straightforward consequence of Lemma \ref{le:A in A}

\begin{lemma}\label{le:D in D}
For any $4\leq p\leq n$,
$$\Alt_{B_{n-p+1}^Y}(A_{B_{n-p+2}^Y}^{c^{(\lambda,\mu)}})=
\left\{\begin{array}{cl}
0&\text{if $1\in\lambda$ or $\mu_i\in 2\IN^*\thickspace\forall i$}\\[1.2ex]
(-1)^{|\lambda|}\cdot
A_{B_{n-p+1}^Y}^{c^{(\lambda\cupone,\mu)}}&\text{otherwise}
\end{array}\right.$$
\end{lemma}
Since $d^\#_{p-1}(f)=-\Alt_{B_{n-p+1}^Y}(f)$ for any $f\in \alg{W_{B_{n-p+2}^Y}}
^\veps$, Lemma \ref{le:D in D} accounts for the third summand in
\eqref{eq:Im D}.
\end{proof}

\subsubsection{}


\begin{proposition}\label{pr:Ker D}
The following holds for any $2\leq p\leq n$,
\begin{equation}\label{eq:Ker D}
\Ker d^\#_p=
\Im d^\#_{p-1}\oplus
\bigoplus_{\substack{\nu\vdash p/2:\\ [\nu]\in 2\IN}}
k A^{c^{(\emptyset,2\nu)}}_{B^Y_{n-p+1}}
\oplus\left\{\begin{array}{cl}
\displaystyle{\bigoplus_{\lambda\in\CO_{p+1}^*}k\wt{B}_\lambda}
&\text{if $p\leq n-1$}\\[1.8 em]
\displaystyle{
\bigoplus_{\substack{0\leq m<n,\\\lambda\in\CO_m^*,\\\mu\in\CO_{n-m}:[\mu]\in 2\IN}}
\negthickspace\negthickspace\negthickspace\negthickspace\negthickspace
kA_{B^Y_1}^{c^{(\lambda,\mu)}}}
&\text{if $p=n$}
\end{array}\right.
\end{equation}
where the second summand only arises if $p$ is even and greater
or equal to $3$, $\CO_q^*=\{\lambda\in\CO_q|\ths 1\notin\lambda\}$
for $q\geq 1$, $\CO_0^*=\{\emptyset\}$ and
$$\wt{B}_\lambda=
\sum_{i=1}^{n-p-1}A_{[i,i+p-1]}^{c^\lambda}+
A_{B^+_{n-p}}^{c^\lambda}+A_{B^-_{n-p}}^{c^\lambda}$$
\end{proposition}
\begin{proof}
For $p\geq 3$ even, the subspace spanned by $A_{B_{n-p+1}^Y}
^{c^{(\emptyset,2\nu)}}$, $\nu\vdash p/2$, lies in $\Ker d^\#_p$ by lemma
\ref{le:D in D} and is in direct sum with $\Im d^\#_{p-1}$ by proposition
\ref{pr:Im D}. This accounts for the second summand in \eqref{eq:Ker D}.
Let now $f\in\Ker d^\#_p$ be such that all components of $f_{B_{n-p+1}^Y}$
of type $c^{(\emptyset,2\nu)}$, $\nu\vdash p/2$, are zero. If $p=n$,
$f$ lies in the span of the elements $A_{B^Y_1}^{c^{(\lambda,\mu)}}$,
$(\lambda,\mu)\in\CO_m\times\CO_{n-m}$ and $0\leq m\leq n$. Since
$d^\#_n=0$ and
$$\Im d^\#_{n-1}=
\bigoplus_{\lambda\in\CO^*_n}
k A_{B_1^Y}^{c^{(\lambda,\emptyset)}}\oplus
\bigoplus_{\substack{1\leq m\leq n,\\\overlineambda\in\CO_m,\\\overline{\mu}\in\CO_{n-m}:
\\1\in\overlineambda,[\overline{\mu}]\in 2\IN}}
k A_{B_1^Y}^{c^{(\overlineambda,\overline{\mu})}}$$
by Proposition \ref{pr:Im D}, \eqref{eq:Ker D} holds for $p=n$. Assume
now that $p\leq n-1$. By lemmas \ref{le:A in D}-\ref{le:D in D}, the
restriction of $d^\#_p$ is injective on the span of the elements $A_{B^Y
_{n-p+1}}^{c^{(\lambda,\mu)}}$ where $(\lambda,\mu)\in\CO_m\times\CO
_{p-m}$, $0\leq m\leq p$, are such that $\lambda$ does not contain $
1$ and $\mu$ is nonempty. The corresponding components of $f$ are
therefore zero. Since $A_{B^Y_{n-p+1}}^{c^{(\lambda,\mu)}}$ lies in
$\Im d^\#_{p-1}$ if $1\in\lambda$ by Lemma \ref{le:D in D},
we may therefore assume that $f_{B^Y_{n-p+1}}$ only has
components along $A_{B^Y_{n-p+1}}^{c^{(\lambda,\emptyset)}}$,
$\lambda\in\CO_p^*$. Working modulo the first summand of $\Im d^\#
_{p-1}$ given by Proposition \ref{pr:Im D}, we may further assume
that $f$ lies in the span of
\begin{equation}\label{eq:rest}
\begin{split}
&\bigoplus_{\lambda\in\CO_{p+1}^*}
\bigoplus_{i=1}^{n-p-1}kA_{[i,i+p-1]}^{c^\lambda}\\
\oplus&\bigoplus_{\lambda\in\CO_{p+1}}\left( 
kA_{B^+_{n-p}}^{c^\lambda}\oplus kA_{B^+_{n-p}}^{c^\lambda}\right)
\oplus\bigoplus_{\overlineambda\in\CO^*_p}
k A_{B^Y_{n-p+1}}^{c^{(\overlineambda,\emptyset)}}
\end{split}
\end{equation}
Let $\lambda\in\CO_{p+1}^*$. It readily follows from \eqref{eq:Ker} and \eqref
{eq:d on B+-} that $\wt{B}_\lambda\in\Ker d^\#_p$. Next, applying \eqref{eq:Ker}
to the $c^\lambda$-components of $f$ along $\alg{W_{[i,i+p-1]}}^\veps$, $i=1,
\ldots,n-p-1$ and $\alg{W_{B_{n-p}^\pm}}^\veps$, we see that these components
are equal to $a_\lambda A_{[i,i+p-1]}^{c^\lambda}$ and $a_\lambda A_
{B_{n-p}^\pm}^{c^\lambda}$ respectively, for some constant $a_\lambda$.
Thus, substracting $a_\lambda\wt{B}_\lambda$ to $f$, we may assume
that all these components are equal to zero and therefore that $f$ lies in
\begin{equation}\label{eq:revisited rest}
\bigoplus_{\substack{\lambda\in\CO_{p+1}:\\1\in\lambda}}\left( 
kA_{B^+_{n-p}}^{c^\lambda}\oplus kA_{B^+_{n-p}}^{c^\lambda}
\oplus
k A_{B^Y_{n-p+1}}^{c^{(\lambda\setminus\{1\},\emptyset)}}\right)
\end{equation}
By \eqref{eq:d on B+-} and Lemma \ref{le:D in D} any solution of $d^\#_p f=0$
with values in \eqref{eq:revisited rest} lies in the subspace $\bigoplus_{\lambda
\in\CO_{p+1}:1\in\lambda}V^\lambda$ defined by \eqref{eq:V lambda} and therefore
in $\Im d^\#_{p-1}$.
\end{proof}

\subsubsection{}

\begin{theorem}\label{th:D}
For any $0\leq p<n$, we have
\begin{equation*}
\begin{split}
\dim HD^p(kW_{D_n})_c=
\begin{cases}
1 & \text{ if } (c=c^{(1^{n-p},2\nu)} \text{ for some }
\nu\vdash p/2 \text{ with }[\nu]\in 2\IN) \\
 & \text{ or } (c=c^{(\lambda\cup\{1^{n-p-1}\},\emptyset)} \text{ with }
\lambda\in\CO_{p+1}^*) \\
0 & \text{otherwise}
\end{cases}
\end{split}
\end{equation*}
and
$$\dim HD^n(kW_{D_n})_c=
\begin{cases}
1 &\text{ if }c=c^{\lambda,\mu} \text{ with }
\lambda\in\CO_m^* \text{ and }
\mu\in\CO_{n-m} \text{ and }[\mu]\in2\IN \\
0 & \text{otherwise.}
\end{cases}
$$
\end{theorem}

\subsection{Stabilisation}

Assume that $W$ is of classical type $\sfX=\sfA$, $\sfB$ or $\sfD$.
Then, using the description of the Dynkin cohomology of $kW$ given in this
section, one readily checks that
for $n\geq m$, the restriction map
on Dynkin diagram cohomology defined in \S \ref{ss:restriction} 
$$\mathrm{Res}_{\sfX_m}^{\sfX_n}:
HD^p(kW_{\sfX_n})\to HD^p(kW_{\sfX_m})$$
is an isomorphism for $p\le m$.

\subsection{Exceptional groups}\label{se:exceptional}

\subsubsection{$\sfG_2$}

The Weyl group of type $\sfG_2=\sfI(6)$ was treated in Section \ref{ss:rank 2}.

\subsubsection{$\sfF_4$}

Let $W$ be the Weyl group of type $\sfF_4$ and label the connected
subdiagrams of $D$ by the subintervals of $[1,4]$ with integral endpoints.
Thus, $W_{[1,3]}$ and $W_{[2,4]}$ are of type $\sfB_3$, $W_{[1,2]}$
and $W_{[3,4]}$ are of type $\sfA_2$, and $W_{[2,3]}$ is of type
$\sfB_2$. It follows from Propositions \ref{pr:A sign} and \ref{pr:B sign}
that
\begin{gather*}
\alg{W_{[1,3]}}^\veps=0
\qquad\qquad
\alg{W_{[2,4]}}^\veps=0\\
\alg{W_{[1,2]}}^\veps\simeq k
\qquad\qquad
\alg{W_{[2,3]}}^\veps\simeq k
\qquad\qquad
\alg{W_{[3,4]}}^\veps\simeq k
\end{gather*}
The differential $d^\#$ of $HC$ is therefore equal to
zero, so that
$$\dim HD^p(kW)=\left\{\begin{array}{cl}
0&\text{if $p=0$}\\
0&\text{if $p=1$}\\
3&\text{if $p=2$}\\
0&\text{if $p=3$}\\
\dim \alg{W}^\veps &\text{if $p=4$}\\
0&\text{if $p\geq 5$}
\end{array}\right.$$

\subsubsection{GAP calculations}

For the groups $\sfH_3$, $\sfH_4$, $\sfE_6$, $\sfE_7$ and $\sfE_8$, the
dimension of the cohomology spaces of the complex $HC$ can be readily
computed with the computer algebra package \cite{GAP}. For each $B$, we
 enumerate conjugacy classes of $W_B$ and select those whose elements
 have a centraliser in the kernel of $\veps$. Then, we compute the matrices
 corresponding to the differentials and determine their rank. We indicate the
 results in the following table, where the columns provide $\dim HD^i(kW)$,
$i=2,\ldots,n$.

\vskip .5cm
\hspace{1.3in}
\begin{tabular}{|c|ccccccc|}
\hline
$i=$        & 2 & 3 & 4 & 5 & 6 & 7 & 8\\
\hline
$\sfE_6$ & 1 & 0 & 2 & 0 & 4 &   &\\
$\sfE_7$ & 1 & 0 & 2 & 0 & 7 & 0 &\\
$\sfE_8$ & 1 & 0 & 2 & 0 & 6 & 1 & 17\\
$\sfF_4$ & 3 & 0 & 5 &   &   &   &\\
$\sfG_2$ & 2 &    &   &    &   &   &\\
$\sfH_3$ & 3 & 0 &   &   &   &   &\\
$\sfH_4$ & 3 & 0 & 16 & &  &   &\\
\hline
\end{tabular}

\section{Affine Coxeter groups}
\label{sec:affine}

The methods of Sections \ref{se:affinecox} and \ref{se:finitepart} provide
a partial computation of the Dynkin cohomology of affine Weyl groups. Let
$D_0$ be a finite crystallographic Dynkin diagram with $n$ vertices and
let $D$ be its completion. Let $V$ be the reflection representation of $W
_0=W_{D_0}$ and let $Q\subset V$ be the coroot lattice, so that $W\simeq
Q\rtimes W_0$.

We have canonical Serre duality isomorphisms
$$\Ext^i_{\alg{Q}}(k,M)\iso \Ext^{n-i}_{\alg{Q}}(M,k)^*\otimes \Lambda^n V^*$$
for any finitely generated $\alg{Q}$--module $M$ and any integer $i$.
Thus, given a finitely generated $\alg{W}$--module $M$ and an integer
$i$, we have isomorphisms
\begin{equation}\label{eq:isoms}
\Ext^i_{\alg{W}}(k,M)\simeq
\Ext^i_{\alg{Q}}(k,M)^{W_0}\simeq 
\left(\Ext^{n-i}_{\alg{Q}}(M,\veps)^*\right)^{W_0}\simeq
\Ext^{n-i}_{\alg{W}}(M,\veps)^*
\end{equation}

\smallskip
Let us describe the conjugacy classes of $W$. Given $v\in W_0$, let
$Q_v \subseteq Q$ be the sublattice given by
$$Q_v=\{x-v(x)\ |\ x\in Q\}$$

\begin{proposition}
Representatives of conjugacy classes of $W$ are given by elements $tv$
where $v$ runs over representatives of conjugacy classes of $W_0$ and
$t$ runs over representatives of $(Q/Q_v)/C_{W_0}(v)$.
\end{proposition}

\begin{proof}
Given $t,t'\in Q$ and $v,v'\in W_0$, 
if $tv$ and $t'v'$ are conjugate in $W$, then $v$ and $v'$ are conjugate in
$W_0$. Let then $g=\tau x$ with $\tau\in Q$ and $x\in W_0$ be such that
$gtvg^{-1}=t'v$. Then, $x\in C_{W_0}(v)$. We have
$gtvg^{-1}=xtx^{-1}\cdot \tau \cdot v\tau^{-1}v^{-1}\cdot v$ and the proposition
follows.
\end{proof}

Let $c$ be a conjugacy class of $W$, let $w\in c$ and let  $\bar{w}$ be the
image of $w$ in $W_0$. The quotient $\overline{C_W(w)}=C_W(w)/Q^{\bar
{w}}$ is a subgroup of $C_{W_0}(\bar{w})$. Denoting as customary the vector
space spanned by the elements of $c$ by $k(c)$, we have
\begin{align*}
\Ext^i_{\alg{W}}(k,\vec{c})
&\simeq \Ext^{n-i}_{\alg{W}}(\vec{c},\veps)^*\\[1.3 ex]
&\simeq \Ext^{n-i}_{\alg{C_W(w)}}(k,\veps)^*\\[1.1 ex]
&\simeq \left(\Ext^{n-i}_{\alg{Q^{\bar{w}}}}(k,k)^*
 \otimes\veps\right)^{\overline{C_W(w)}}\\[1.1 ex]
&\simeq \left(\Lambda^{n-i}(V^{\bar{w}})
 \otimes\veps\right)^{\overline{C_W(w)}}.
\end{align*}
where the first isomorphism uses \eqref{eq:isoms}, the second one
Frobenius reciprocity and the last one the fact that, given a
finitely--generated free abelian group $L$, we have the Koszul
isomorphism
$H^*(L,k)\simeq \Lambda^*(L^*\otimes_\IZ k)$.

\begin{theorem}
If the elements in $c$ have infinite order, then
$$HD^i_{c}\simeq \left(\Lambda^{n+1-i}(V^{\bar{w}})
 \otimes\veps\right)^{\overline{C_W(w)}}.$$
\end{theorem}

\begin{proof}
The assumption on $c$ shows that $c\cap W_B=\emptyset$ for any
proper subset $B$ of $D$. Thus, the subcomplex $CD_{c}$ defined
in Section \ref{ss:conjugacy} is concentrated in degree $n$. The
theorem now follows from \S \ref{se:affinecox} and the isomorphisms
above.
\end{proof}

\begin{remark}
If the elements of $c$ have finite order, then $c$ has a non empty
intersection with $W_B$ for some proper connected subdiagram
$B$ of $D$ \cite{Hu}. In that case, there is a distinguished triangle
$$\bigoplus_i \left(\Lambda^{n-i+1}(V^{\bar{w}})
 \otimes\veps\right)^{\overline{C_W(w)}}[-i]\to CD_c \to
{HC_f}_{c}\rightsquigarrow$$
where ${HC_f}_{c}$ is the subcomplex of $HC_f$ given by
${HC_f}_{c}^i=\bigoplus_B \vec{c\cap W_B}^\veps$ and $B$
runs over the proper subdiagrams of $D$ of size $i$.
\end{remark}

\end{document}